\documentclass[final,onefignum,onetabnum]{siamsials251208}

\usepackage{lipsum}
\usepackage{amsfonts}
\usepackage{graphicx}
\usepackage{epstopdf}
\usepackage{algorithmic}
\ifpdf
  \DeclareGraphicsExtensions{.eps,.pdf,.png,.jpg}
\else
  \DeclareGraphicsExtensions{.eps}
\fi

\usepackage{enumitem}
\setlist[enumerate]{leftmargin=.5in}
\setlist[itemize]{leftmargin=.5in}

\newsiamremark{remark}{Remark}
\newsiamremark{hypothesis}{Hypothesis}
\crefname{hypothesis}{Hypothesis}{Hypotheses}
\newsiamthm{claim}{Claim}
\newsiamremark{fact}{Fact}
\crefname{fact}{Fact}{Facts}

\headers{Beyond the Critical Depth}{M. A. Neto, P. A. Marquet, M. A. Freilich, L. Martí, and N. Sanchez-Pi}

\title{Beyond the Critical Depth: The Metabolic and Physical Drivers of Phytoplankton Persistence in a Changing Ocean%
\thanks{Submitted to SIAM Journal on Life Sciences
\funding{This work was supported by the Franco-Chilean Binational Center of Artificial Intelligence, ANID Strengthening R\&D capabilities Program CTI230007 Inria Chile, and Inria Challenge OcéanIA (desc. num 14500).}}}

\author{Matías A. Neto\thanks{Inria Chile Research Center, Las Condes, RM, Chile 
  (\email{matias.neto@inria.cl}, \email{luis.marti@inria.cl}, \email{nayat.sanchez-pi@inria.cl}).}
\and Pablo A. Marquet\thanks{Facultad de Ciencias Biol\'ogicas, Pontificia Universidad Católica de Chile, Santiago; Santa Fe Institute, NM, USA; Centro de Modelamiento Matemático (CMM), Chile; Centro de Cambio Global UC, Chile; Instituto de Sistemas Complejos de Valpara\'iso (ISCV), Chile (\email{pmarquet@uc.cl}).}
\and Mara A. Freilich\thanks{Division of Applied Math and Department of Earth, Environmental, and Planetary Sciences, Brown University, Providence, RI, USA (\email{mara\_freilich@brown.edu}).}
\and Luis Martí\footnotemark[2]
\and Nayat Sanchez-Pi\footnotemark[2]}

\usepackage{amsopn}

\ifpdf
\hypersetup{
  pdftitle={Beyond the Critical Depth: The Metabolic and Physical Drivers of Phytoplankton Persistence in a Changing Ocean},
  pdfauthor={M. A. Neto, P. A. Marquet, M. A. Freilich, L. Martí, N. Sanchez-Pi}
}
\fi

\usepackage{xspace}

\usepackage[version=4]{mhchem}
\usepackage{graphicx}%
\usepackage{subcaption}
\usepackage{multirow}%
\usepackage{xcolor}%
\usepackage{textcomp}%
\usepackage{manyfoot}%
\usepackage{booktabs}%
\usepackage{listings}%
\usepackage{rotating}
\usepackage{comment}
\usepackage{siunitx}
\usepackage{multiobjective}
\usepackage{xfrac}
\usepackage{makecell}
\usepackage{cite}
\usepackage[showdeletions]{color-edits}
\addauthor[Matias]{mn}{orange}
\addauthor[Pablo]{pm}{magenta}
\addauthor[Nayat]{nspi}{BrickRed}
\addauthor[Luis]{lm}{cyan}
\addauthor[Mara]{mf}{ForestGreen}

\usepackage[T1]{fontenc} 
\newcommand{\etal}{\emph{et\,al.}\xspace}
\newcommand{\GFDLStableViabilityArea}{318.89}
\newcommand{\GFDLStableViabilityPerc}{91.7}
\newcommand{\GFDLStableRestrictionArea}{10.18}
\newcommand{\GFDLStableRestrictionPerc}{2.9}
\newcommand{\GFDLHabitatExpansionArea}{8.55}
\newcommand{\GFDLHabitatExpansionPerc}{2.5}
\newcommand{\GFDLHabitatContractionArea}{0.40}
\newcommand{\GFDLHabitatContractionPerc}{0.1}
\newcommand{\GFDLIceFreeViabilityArea}{2.12}
\newcommand{\GFDLIceFreeViabilityPerc}{0.6}
\newcommand{\GFDLIceFreeRestrictionArea}{7.52}
\newcommand{\GFDLIceFreeRestrictionPerc}{2.2}
\newcommand{\GFDLHyperThreshold}{0.80}

\newcommand{\GFDLHyperAboveThresholdPerc}{83.1}
\newcommand{\GFDLHyperMeanRsq}{0.859}

\begin{document}
\maketitle

%%==================================%%
\begin{abstract} While the classical Critical Depth Hypothesis (CDH) effectively explains the onset of blooms as transient instabilities, it does not fully capture the seasonal decoupling of biological rates and the long-term persistence of phytoplankton communities in fluctuating thermal environments. To address these limitations, we introduce a parsimonious framework that leverages the theory of non-autonomous dynamical systems to diagnose the stability of phytoplankton communities throughout the entire annual cycle. By linearizing the dynamics around the extinction equilibrium, we identify the invasion growth rate -- formally the Floquet exponent-- and derive the critical nutrient requirement ($\gamma_\text{crit}$) as a bifurcation point for uniform persistence. Using end-of-the-century projections from the GFDL-ESM4 model under a high-emission scenario (SSP5-8.5), we identify a global regime shift characterized by a widespread expansion of metabolic-driven regimes, which increasingly displace regions where stability was historically governed by physical mixing.
\end{abstract}
\begin{relevance}
Quantitative analysis of system stability challenges CDH by demonstrating that metabolic constraints increasingly modulates phytoplankton persistence in a changing ocean. Our results, based on high-emission projections, reveal a profound physical-biological decoupling at the poles: while warming reduces the critical nutrient requirement ($\gamma_\text{crit}$) facilitating persistence in previously marginal waters, this metabolic expansion is offset at poles. A 1:4 ratio between newly viable niches and ice-free deserts suggests that cryospheric retreat does not guarantee a proportional expansion of life. In addition, we identify the North Atlantic Subpolar Gyre as a ``metabolic refuge'' where mixing dynamics still anchor the ecosystem against global thermalization. By providing a ``radiography'' of the future ocean's complexity, this methodology offers a mechanistic basis to deconstruct how the dynamic balance between environmental energy and metabolic demands may determine the functional integrity of the marine biosphere under extreme anthropogenic forcing.
\end{relevance}
\begin{mathcontent}
The temperature dependence of biological rates is modeled using a thermodynamic equation, coupling population dynamics with seasonal variations in mixed layer depth and temperature. Given the non-autonomous nature of the system under annual forcing, we characterize the stability of the extinction equilibrium through its associated invasion growth rate. This rate is analytically derived as the Floquet exponent $\lambda_P$, which provides a rigorous condition for uniform persistence (Theorem \ref{th:persistence}). The numerical analysis of this exponent, projected onto a global scale, quantifies the relative influence of environmental drivers on the stability threshold $\gamma_\text{crit}$. This allows for the definition of the thermal dominance index ($\mathcal{D}_T$), a metric that identifies the geographic transition from mixing-driven to metabolic-driven ecological control.
\end{mathcontent}

% REQUIRED
\begin{keywords}
Phytoplankton persistence, Non-autonomous dynamics, Critical nutrient requirement, Metabolic dominance, Climate change
\end{keywords}

% REQUIRED
\begin{MSCcodes} 
37N25, 37B55, 92D40, 37C75, 37N10
\end{MSCcodes}



% \pagestyle{fancy}
% \fancyhead{} % clear all header fields
% \renewcommand{\headrulewidth}{0pt} % no line in header area
% \fancyfoot{} % clear all footer fields
% \fancyfoot[LE,RO]{\thepage}           % page number in "outer" position of footer line
% \fancyfoot[RE,LO]{\textcolor{red}{\footnotesize Work in progress, do not disseminate. Comments welcomed!}}

\section{Introduction}\label{sec1}

Marine microbial communities play a central role in global biogeochemistry, regulating atmospheric and oceanic chemistry through photosynthesis and carbon sequestration \cite{falkowski2008microbial,fuhrman2015marine}. While phytoplankton represent only about $1\%$ of global plant biomass, they account for more than half of global photosynthetic activity, driving the biological carbon pump by exporting atmospheric \ce{CO2} into the deep ocean \cite{henson2011, litchman2012phytoplankton}. Traditionally, the coupling between biological productivity and the physical environment has been explained by theories that focus on how light limitation and mixed layer depth (MLD) trigger biomass accumulations known as blooms, such as the Critical Depth Hypothesis (CDH) \cite{sverdrup} or Dilution Recoupling Hypothesis \cite{behrenfield2010}. However, as the global ocean undergoes rapid and heterogeneous changes -- ranging from accelerated warming in polar regions \cite{meredith2019polar} to anomalous cooling in subpolar gyres \cite{rahmstorf_2015} -- it is increasingly important to examine the impacts of temperature on the long-term persistence of these bloom-forming ecosystems.

Given the critical role of phytoplankton communities, deriving rigorous conditions for their stability under shifting climate scenarios has become a pressing scientific priority. Current Earth System Models (ESMs) provide essential projections of ecosystem responses, yet their internal complexity can obscure the mechanistic link between physical forcing and biological processes \cite{anderson_2005}. Given the bloom-forming nature of marine phytoplankton communities, it is critical to understand how the nature of periodic dynamics may change in response to physical forcing and to develop general theories that can account for periodicity.  The general theory of periodic systems provides a robust background \cite{zhao,wang2008threshold} -- and has been explored in specific ecological models \cite{MingChen2017,klausmeier_floquet_2008} -- a formal analytical bridge to global marine ecosystems remains absent.

We introduce a parsimonious framework that treats the marine ecosystem as a non-autonomous dynamical system. Rather than focusing on instantaneous balances, we evaluate the annual stability of the ecosystem by calculating the invasion growth rate -- formally derived as the Floquet exponent. This approach, allows us to define a structural persistence threshold, $\gamma_\text{crit}$, which represents the ``metabolic tax'' imposed by the environment: the minimum nutrient level required for a population to recover from its seasonal minimum and persist. By linearizing the dynamics around the extinction equilibrium, we obtain a single metric that accounts for the integrated periodic effects of temperature-dependent metabolism and physical mixing constraints.

The remainder of this paper is organized as follows. In Section \ref{sec:model}, we formulate an NPZ-T model that incorporates thermodynamics theory \cite{arroyo} into seasonal physics. In Section \ref{sec:threshold}, we derive the analytical stability criteria for the extinction equilibrium and define the threshold for phytoplankton persistence. Finally, in Section \ref{sec:global_analysis}, we apply this framework to global ESM projections to map persistence boundaries through the 21st century. By examining these boundaries under a high-emission scenario (SSP5-8.5), we test the limits of ecosystem resilience and the emerging dominance of metabolic constraints.

The main advance of this work is twofold: from a diagnostic perspective, it provides a rigorous analytical bridge to quantify ecosystem stability in complex climate models; biologically, it reveals a shift in the controls on marine ecosystems where the fate of phytoplankton is increasingly dictated by metabolic demands rather than physical mixing. This framework offers a mechanistic lens to evaluate the resilience of the marine biosphere in a changing world.

\section{NPZ-T: Temperature-Dependent NPZ Model with MLD}\label{sec:model}

We define the NPZ-T model as an extension of the classical NPZ models \cite{franks,evans}, where the physical dynamics of MLD and the temperature-dependent biological rates are explicitly coupled with the biological interactions of nutrients, phytoplankton and zooplankton. 

Let $N(t)$, $P(t)$, and $Z(t)$ denote the nutrient, phytoplankton, and zooplankton concentrations, respectively. The dynamics are described as
\begin{equation}\label{eqn:NPZ-T-full}
\begin{aligned}
\frac{dN}{dt} &= - f(t,I,H) g(t,N,T) P + (1-\alpha) h(t,P,T) Z + m_P(t,T) P + m_Z(t,T) Z^2 \\
&+ s_+(t) (P+Z),\\
\frac{dP}{dt} &= f(t,I,H) g(t,N,T) P - h(t,P,T) Z - m_P(t,T) P - s_+(t) P\,,\\
\frac{dZ}{dt} &= \alpha h(t,P,T) Z - m_Z(t,T) Z^2 - s_+(t) Z\,,
\end{aligned}
\end{equation}
where $I$ is the surface irradiance, $H$ is the MLD, $T$ is the Sea Surface Temperature (SST), $f(\cdot)$ represents light absorption, $g(\cdot)$ captures the nutrient uptake, $h(\cdot)$ expresses the  zooplankton grazing. Parameters $m_P(\cdot)$ and $m_Z(\cdot)$ denote mortality rates, $s_+$ accounts for dilution from entrainment of water below the mixed layer, and $\alpha$ denotes the zooplankton assimilation efficiency. This is a bulk mixed layer model, which models average dynamics within the mixed layer and neglects lateral gradients \cite{evans,Mara}.

\subsection{Physical Forcing and Seasonal Dynamics} 
The average light limitation within the MLD and the dilution effects due to deepening are modeled as
\begin{align}
f(t,I,H) & = \frac{1}{K_d H(t)}\left(1 - e^{-K_d H(t)}\right) \frac{I(t)}{I(t)+I_0}\,,\\
s_+(t) & = 
\begin{cases}
\frac{1}{H(t)} \frac{d H(t)}{dt}, & \text{if }\frac{d H(t)}{dt} \ge 0,\\
0, & \text{otherwise.}
\end{cases},
\end{align}

\noindent where $K_d$ is the light attenuation coefficient, $I_0$ is the half-saturation constant and $s_+(t)$ captures the dilution of biomass as the mixed layer deepens. Physically, $s_+$ represents entrainment of nutrient-rich deep water: as the mixed layer expands, it recycles the diluted organic nutrient back into the inorganic pool $N$. This mechanism reflects the seasonal nutrient replenishment while ensuring mass conservation within the system.

Following the analytical approach of Freilich~\etal \cite{Mara}, we approximate the MLD, $H(t)$, and temperature $T(t)$ using periodic functions that capture the characteristic minima and maxima of the seasonal cycle. These variables are derived by fitting the periodic functions to seasonal climatologies from GFDL-ESM4 model \cite{gfdl-esm4} at each grid point. Specifically, for $T(t)$, the fit captures the amplitude, phase, and mean values, whereas for $H(t)$, the formulation follows \cite{Mara} by utilizing only the phase and the maximum annual depth, effectively anchoring the seasonal expansion to its winter peak. To ensure the reliability of the seasonal cycle fits and focus on ice-free pelagic dynamics, we exclude grid points where the minimum annual temperature falls below $\qty{-1.8}{\unit{\degreeCelsius}}$. In these regions, cryospheric processes and extreme light limitation -- not explicitly modeled here -- significantly alter the biological response. The seasonal surface irradiance, $I(t)$, follows the formulation of Brock \cite{Brock1981} to accurately reproduce the latitudinal modulation of the geometrical solar cycle.

\subsection{Nutrient Uptake, Grazing, Mortality and Temperature Dependence}

Nutrient limitation is modeled using Monod kinetics, while zooplankton grazing follows a Holling Type III functional response, consistent with the dilution-recoupling mechanism \cite{Mara,behrenfield2010}
\begin{align}
g(N,T) & = \mu_P(T) \frac{N}{N + N_0}\,,\\
h(P,T) & = \mu_Z(T) \frac{P^2}{P^2 + P_0^2}\,,
\end{align}
where $\mu_P(T)$ and $\mu_Z(T)$ denote the temperature-dependent maximum growth and grazing rates, respectively, and $N_0$ and $P_0$ are half-saturation constants. Phytoplankton and zooplankton mortality rates, $m_P(T)$ and $m_Z(T)$, are also temperature-dependent.

To explicitly capture the physiological non-linearities and thermal denaturation we adopted the general thermodynamic framework proposed by Arroyo \etal \cite{arroyo}. Unlike standard exponential formulations (e.g., $Q_{10}$ or Arrhenius) used in previous NPZ models \cite{ROY201216, Gregg2008}, this framework constitutes a generalization of the classical Eyring-Evans-Polanyi (EEP) theory that can produce concave and convex Thermal Performance Curves (TPCs). We express the metabolic rate function normalized to a reference temperature $T_0=\qty{293.15}{\kelvin}$ (\qty{20}{\degreeCelsius}) as
\begin{equation}\label{eqn:exp-eep}
k(T)=k(T_0)\exp\left[\frac{\Delta C}{R}\ln{\left(\frac{T}{T_0}\right)}+\frac{\Delta H}{R}\left(\frac{1}{T_0}-\frac{1}{T}\right)\right]\,,
\end{equation}
where $k(T_0)$ is the metabolic rate at the reference temperature, $R$ is the universal gas constant. The thermodynamic parameters represent the change in heat capacity of activation ($\Delta C$), the enthalpy of activation ($\Delta H$). This is to ensure numerical stability in our computational simulations and to provide clear biological meaning to the baseline rates.

Although rooted in chemical kinetics, this formulation scales robustly to ecological level, describing the thermal response of organism from microorganism to metazoans. To ensure ecological consistency, both phytoplankton and zooplankton are modeled as concave ectotherms following the thermodynamic constraints of Arroyo \etal \cite{arroyo}.

\subsection{Phytoplankton Parameters}

\begin{figure}[tb]
    \centering
    \includegraphics[width=0.9\textwidth]{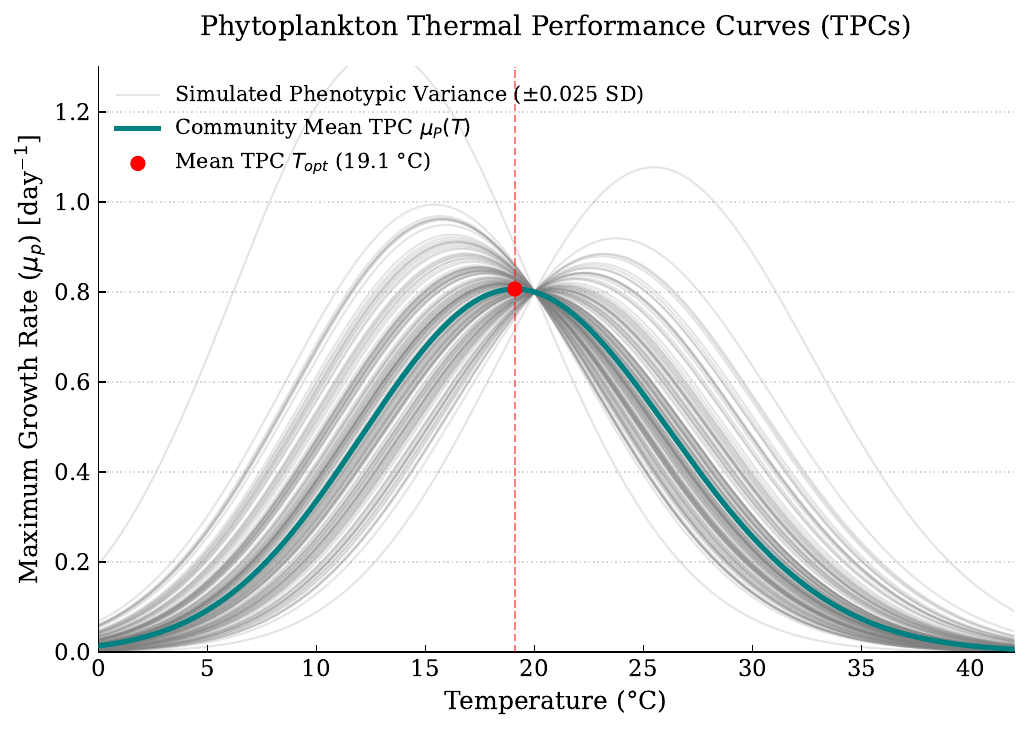}
    \caption{Thermal Performance Curve (TPC) for the maximum phytoplankton growth rate (\(\mu_P\)) derived from the the thermodynamic model \ref{eqn:exp-eep}. The solid teal line represents the community mean TPC used for the global simulations. The background gray lines show a simulated phenotypic ensemble (generated by randomly sampling \(\Delta C\) and \(\Delta H\) with a variance scaled to \(0.025\) SD of the global empirical data), representing the thermal flexibility of local communities.}
    \label{fig:gp_temperature_dependence}
\end{figure}

We parametrize a generic phytoplankton functional type with a maximum growth rate of $\mu_P=\qty{0.8}{\per\day}$ at \qty{20}{\degreeCelsius} \cite{Eppley1972, bissinger2008predicting}. The TPC follows a concave trajectory characterized by the following thermodynamic constraints (mean $\pm$ standard deviation)
\begin{equation}
    \Delta C_P = -14418.93 \pm 2879.40 \,\mathrm{J}, \quad \Delta H_P = 4213943.43 \pm 838418.23 \,\mathrm{J}.
\end{equation}
While individual phytoplankton species often display highly asymmetric TPCs with an abrupt physiological decline past their thermal optimum, the immense phenotypic diversity across global oceans results in large standard deviations in these parameters. Averaging these traits to represent a global, generic phytoplankton community naturally smooths out individual asymmetries, yielding a more symmetric Community mean TPC (Figure \ref{fig:gp_temperature_dependence}). Mortality is set to $m_P=\qty{0.024}{\per\day}$ at $20^\circ$C \cite{Baker2021}.

\subsection{Zooplankton Parameters}

Reference rates at \qty{20}{\degreeCelsius} of a zooplankton generic type are set at $\mu_Z=\qty{4}{\per\day}$ and $m_Z=\qty{1.8}{\per\day}$ \cite{Pulsifer2021, Mara}. The thermodynamic curvature is defined by the same thermodynamic constraints of the phytoplankton.

The resulting temperature-dependent rates for both populations are computed dynamically at each time step based on the local temperature profile. Table \ref{tab:model_parameters} summarizes the variables and constants used in the model.

\begin{table}[tb]
    \centering
    \caption{Definition of model variables and parameters.}
    \label{tab:model_parameters}
    \begin{tabular}{lll}
        \toprule
        \textbf{Symbol} & \textbf{Description} & \textbf{Units} \\
        \midrule
        $H(t)$ & Mixed layer depth & \unit{\metre} \\
        $T(t)$ & Sea Surface Temperature & \unit{\degreeCelsius} \\
        $K_d$ & Light attenuation coefficient & \unit{\per\metre} \\
        $I(t)$ & Surface irradiance & \unit{\micro\mole\per\metre\squared\per\second} \\
        $N_0, P_0$ & Half-saturation constants & \unit{\mmol\per\metre\tothe{3}} \\
        $\mu_P(T), \mu_Z(T)$ & Max. growth/grazing rates (daily) & \unit{\per\day} \\
        $m_P(T), m_Z(T)$ & Mortality rates (daily) & \unit{\per\day} \\
        $\alpha$ & Zooplankton assimilation efficiency & -- \\
        \bottomrule
    \end{tabular}
\end{table}

\section{Threshold Dynamics and Global Persistence Analysis} \label{sec:threshold}

The long-term fate of the marine ecosystem in our model is determined by the stability of the extinction equilibrium. While this approach generalizes the critical depth criteria of Sverdrup \cite{sverdrup}, it is formally rooted in the modern theory of coexistence and invasion analysis \cite{Chesson2000}. Our framework is analogous to the derivation of the Basic Reproduction Number ($R_0$) in epidemiological models for periodic environments \cite{wang2008threshold}; we investigate the per-capita growth rate of a population when rare, identifying the conditions under which a small plankton population can invade and persist in an environment where nutrients are at their maximum availability.

The link between local instability and global persistence is fundamental to our analysis. Although the full non-linear dynamics of the NPZ-T model can be complex, the fundamental condition for any biological activity to occur is the instability of the extinction equilibrium. In a seasonally forced environment, we characterize this stability through the following regime classification
\begin{itemize}
    \item \emph{Global extinction:} If the extinction equilibrium is stable, the population cannot grow even under optimal nutrient availability, leading biomass to vanish for any initial condition.
    \item \emph{Uniform persistence:} If the extinction equilibrium is unstable -- indicated by a positive invasion growth rate -- the population is repelled from the zero boundary. As demonstrated in Section \ref{sec:proofs}, this instability is sufficient to guarantee that the system remains bounded away from zero, ensuring ecological long-term ecological persistence.
\end{itemize}

\subsection{Analytical Determination of the Invasion Rate}

Despite the seasonal variability of the environment, the system possesses a unique invariant boundary solution: the unexploited equilibrium $E_0=(C_0,0,0)$, where $C_0=N(0)+P(0)+Z(0)$ is the total nutrient inventory in the absence of biomass. To determine the conditions for invasion, we consider a small perturbation $(P,Z)$ around $E_0$.

A structural property of our framework is that the coupling between phytoplankton and zooplankton vanishes upon linearization at the origin. As long as the functional responses are of class $\mathcal{C}^1$ at $E_0$, the predator pressure becomes negligible when prey is scarce. Thus, the linearized dynamics follow a diagonal periodic system
\begin{equation}\label{local-system}
    \frac{d}{dt} \begin{pmatrix} P \\ Z \end{pmatrix} = J(t) \begin{pmatrix} P \\ Z \end{pmatrix}, \quad
    J(t) = 
    \begin{bmatrix}
      f(t,I,H)g(t,C_0,T) - m_P(t,T) - s_+(t) & 0 \\
      0 & - s_+(t)
    \end{bmatrix}\,.
\end{equation}
In this diagonal structure, the stability of $E_0$ is determined by Floquet theory \cite{Hale}. The $E_0$ is locally asymptotically stable if and only if all Floquet exponents are negative. Conversely, the equilibrium is unstable if at least one exponent is positive, providing the condition for population invasion. For a detailed mathematical foundations on the concepts underpinning Floquet theory and persistence, we refer the reader to Appendix \ref{sec:foundations}. 

Since the linearized system \eqref{local-system} is diagonal, its fundamental matrix solution $\Phi(t)$ can be obtained through direct integration. The resulting monodromy matrix $C = \Phi(T_p)$ is given explicitly by
\begin{equation*}
    C = 
    \begin{bmatrix}
    \exp\left(\int_0^{T_p} [f(s,I,H)g(s,C_0,T)-m_P(s,T)-s_+(s)] \, ds\right) & 0\\
    0 & \exp\left(-\int_0^{T_p} s_+(s) \, ds\right)
    \end{bmatrix}
\end{equation*}
where the diagonal elements correspond to the Floquet multipliers for phytoplankton ($\rho_P$) and zooplankton ($\rho_Z$), respectively. As the loss rate $s_+(t)$ is strictly positive, the integral for the second multiplier is always negative, ensuring that $|\rho_Z| < 1$. This confirms that zooplankton perturbations always decay in the absence of prey. Therefore, the stability of the extinction equilibrium depends exclusively on the phytoplankton multiplier $\rho_P$ (the spectral radius of the monodromy matrix is defined as $r(C)=\rho_P$). This allows us to define the stability threshold through the invasion growth rate $\lambda_P$:
\begin{equation}
  \lambda_P = \frac{1}{T_p} \ln(\rho_P) = \frac{1}{T_p}\int_0^{T_p} \left( f(s,I,H)g(s,C_0,T) - m_P(s,T) - s_+(s) \right) \, ds\,.  
\end{equation}

\subsection{Global Thresholds for Extinction and Persistence}

It is now critical to bridge the gap between the local stability analysis and the global asymptotic behavior of the full NPZ-T model. This relationship is formalized in the following two theorems.

\begin{theorem}[Global Extinction Condition]\label{th:extinction}
Consider the NPZ-T system \eqref{eqn:NPZ-T-full} with any positive initial condition $(N(0),P(0),Z(0))$ such that $N+P+Z=C_0$. If $\lambda_P<0$, then $E_0$ is globally asymptotically stable, \emph{i.\,e.},
\begin{equation}
  \lim_{t\to \infty} P(t)= 0 \text{ and } \lim_{t\to \infty} Z(t)= 0\,.  
\end{equation}

\end{theorem}

This theorem establishes a fundamental physical barrier: If $\lambda_P<0$, physiological losses and physical dilution outpace the maximum potential growth even under saturating nutrient conditions. Hence, if a population cannot sustain itself even when resource limitation is removed, it will inevitably collapse under the nutrient-limited conditions of the global ocean. 

\begin{theorem}[Global Persistence Condition]\label{th:persistence}
Consider the NPZ-T system \eqref{eqn:NPZ-T-full} with any positive initial condition $(N(0),P(0),Z(0))$ such that $N+P+Z=C_0$. If $\lambda_P>0$, then the system is uniformly persistent within the invariant set defined by $C_0$. Specifically, there $\exists~\eta>0$, such that all solutions satisfy
\begin{equation}
    \liminf_{t\to \infty}{P(t)}\geq \eta\,.
\end{equation}

\end{theorem}

The existence of $\eta>0$ guaranties ecological resilience; the physical and metabolic balance remains positive at the unexploited state. The phytoplankton population will not only survive but will recover to stable densities after environmental perturbations significantly reduce the concentration.

\begin{remark}
    The correspondence between the local instability of $E_0$ (whenever $\lambda_P>0$) an the global property of uniform persistence is consistent with the general threshold theory of periodic compartmental models \cite{wang2008threshold}. Specifically, \cite{chen} demonstrate that for seasonal systems, the condition $\rho(\Phi)>1$ -- where $\rho(\Phi)$ is the spectral radius of the evolution operator -- is equivalent to an ecological reproductive index $R_0>1$, which ensures the transition from extinction to a stable periodic bloom. In our framework, the structural decoupling at the extinction state -- where the nutrient remains at a constant baseline $C_0$ -- ensures that $\lambda_P$ is determined solely by the annual reproductive success against a stable resource background. This formally anchors the invasion growth rate $\lambda_P$ in the theory of periodic dynamical systems, providing a rigorous generalization of the classical Critical Depth criteria to non-autonomous nutrient-limited environments.
\end{remark}

\subsection{Beyond the Critical Depth Hypothesis (CDH): A Metabolic Extension}

The classical CDH \cite{sverdrup} establishes that bloom initiation occurs when the MLD is shallower than a physical horizon where integrated photosynthesis exceeds respiration ($Z_c$). By incorporating thermally modulated metabolic rates, our framework extends this logic, revealing that ecological viability is not solely determined by this light-mixing balance, but also by a non-linear metabolic boundary.

This introduces a complex eco-physiological tension. On the metabolic side, while warming can accelerate growth kinetics, it simultaneously imposes higher metabolic costs and thermodynamic constraints. On the physical side, fluctuations in the MLD dictate light availability while also driving physical dilution—which, rather than acting solely as a loss, can actually promote net growth by decoupling phytoplankton from grazing pressure. Therefore, $\lambda_P>0$ defines the threshold where these competing effects yield a net positive balance. Under climate change, this implies that even in physically optimal conditions, the acceleration of metabolic demands can override photosynthetic gains. This shift effectively redefines the ocean's stability not as a purely physical horizon, but as a dynamical competition between environmental energy and metabolic demand.

\section{Global Boundaries of Persistence}\label{sec:global_analysis}

The theoretical conditions derived in Theorems \ref{th:extinction} and \ref{th:persistence} provide a binary criterion for survival based on the sign of the invasion exponent $\lambda_P$. However, evaluating this exponent globally presents a challenge: while physical variables are well-constrained, the specific nutrient inventory $C_0$ remains highly uncertain in global projections \cite{kwiatkowski, seferian2020tracking}. 

To overcome this limitation, we adopt an inverse persistence framework. Instead of pre-supposing a nutrient load, we calculate the minimum nutrient uptake efficiency required to prevent extinction. This allows us to quantify the thermodynamic cost-of-living imposed by the environment independently of chemical uncertainties. 

\subsection{The Critical Nutrient Requirement ($\gamma_\text{crit}$)}
Recall that the invasion exponent $\lambda_P$ represents the net annual growth rate. By separating gains from the losses, the persistence condition $\lambda_P > 0$ can be rewritten as
\begin{equation*}
    g(C_0) \underbrace{\int_0^{T_p} \mu_P(t,T)f(t,I,H) \, dt}_{\mathcal{G}(x)} > \underbrace{\int_0^{T_p} (m_P(t,T) + s_+(t)) \, dt}_{\mathcal{L}(x)},
\end{equation*}

We define the Critical Nutrient Requirement, $\gamma_\text{crit}(x)$ as the nutrient saturation level $g(C^*_0)$ necessary to reach the break-even point ($\lambda_P=0$): 

\begin{equation}\label{eq:gamma_crit}
    \gamma_\text{crit}(x) = \frac{\mathcal{L}(x)}{\mathcal{G}(x)}\,.
\end{equation}

Because both $\mathcal{L}$ and $\mathcal{G}$ are driven by temperature-dependent rates, $\gamma_\text{crit}$ quantifies the ``metabolic tax'' imposed by the environment. If $\gamma_\text{crit}(x)=0.5$, the population must operate at $50\%$ of its maximum uptake capacity just to survive, any further limitation would lead to local extinction. Hence, persistence is only possible if the local total inventory allows for a saturation level that exceeds the threshold $g(C_0)>\gamma_\text{crit}$.

Using monthly physical output from the GFDL-ESM4 model \cite{gfdl-esm4}, we classify the global ocean into three dynamic regimes

\begin{enumerate}
    \item \emph{Robust Zones ($\gamma_\text{crit} \leq 0.5$).}  Regions where potential growth far exceeds losses ($\mathcal{L} \ll \mathcal{G}$). These areas are structurally resilient, allowing for self-sustaining populations even under suboptimal resource conditions. 
    \item \emph{Marginal Zones ($0.5 < \gamma_\text{crit} \leq 1.0$).} Areas where local persistence is possible but sensitive to environmental fluctuations. These represent potential hot-spots of vulnerability where minor shifts in temperature or mixing can compromise population stability.
    \item \emph{Restrictive Zones ($\gamma_\text{crit} > 1.0$).} Regions where environmental losses exceed the maximum theoretical gain ($\mathcal{L} > \mathcal{G}$), precluding local persistence ($\lambda_P < 0$, Theorem \ref{th:extinction}). These ``metabolic deserts'' identify areas where 1D vertical dynamics alone cannot sustain a population. The presence of biomass in such regions likely points to sub-grid subsidies not captured by this framework, such as lateral advection of nutrients and biomass or specialized physiological adaptations to extreme limitation.
\end{enumerate}

\begin{figure}[tb]
    \centering
    \includegraphics[width=\textwidth]{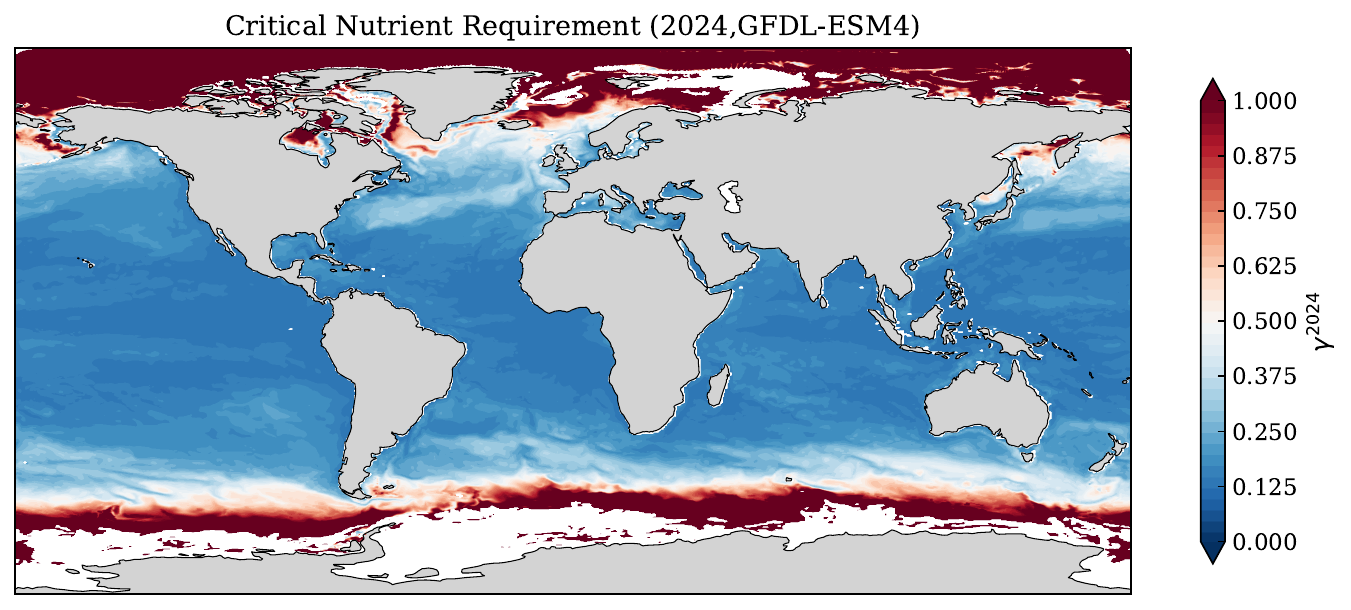}
    \caption{Global Map of the Critical Nutrient Requirement ($\gamma_\text{crit}$) derived from GFDL-ESM4 model data (2024). The color scale quantifies dynamical regime: Blue zones (Robust, $\gamma_\text{crit} \leq 0.5$); Blue to red zones (Marginal, $0.5 < \gamma_\text{crit} \leq 1.0$); Red zones (Restrictive, $\gamma_\text{crit} > 1$)}
    \label{fig:gamma_crit}
\end{figure}

Figure~\ref{fig:gamma_crit} reveals a clear spatial gradient in the ``cost of persistence''. In the subtropics, low $\gamma_{\text{crit}}$ values indicate a low-cost metabolic regime, where ecosystems can persist even under oligotrophic conditions. Conversely, high-latitude regions operate near the theoretical limit ($\gamma_{\text{crit}} \to 1$), meaning they require near-maximal nutrient saturation to offset high physical and metabolic losses. This high baseline cost makes polar and subpolar regions far more sensitive to climatic shifts: because they are already operating near the threshold, any warming-driven increase in metabolic demand or change in mixing can rapidly push the system toward local instability.

\subsection{Metabolic and Physical Drivers of $\gamma_\text{crit}$}

To determine the mechanism governing the shifts in phytoplankton persistence, we decompose the variability of $\gamma_\text{crit}$ into its physical and metabolic components. This is achieved by quantifying the relative contribution of thermal-metabolic effects versus physical mixing constraints through a local sensitivity analysis. 

Let $y$ denote a specific year (e.g., $2024$ or $2100$) and $\mathcal{P}_{\text{ref}}$ the historical reference period (1960--2000). We define the thermal and the mixing anomalies as
\begin{align}
    \Delta T(x)&= \bar{T}_{y}(x) -\langle \bar{T}(x)\rangle_{\text{ref}},\\
    \Delta H(x)&= H^{\max}_{y}(x)-\langle H^{\max}(x)\rangle_{\text{ref}}.
\end{align}
Here, $\bar{T}_{y}$ is the annual mean temperature of the year $y$, $\langle \bar{T}(x)\rangle_{\text{ref}}$ is the climatological mean over the reference period. Similarly, $H^{\max}_{y}$ is the maximum MLD in the year $y$, and $\left\langle H^{\max}(x)\right\rangle_{\text{ref}}$ is the historical mean of the annual maximums. We use $H^{\max}$ as the representative mixing anomaly because it determines the annual nutrient ventilation and the maximum seasonal dilution stress.

To translate these environmental shifts into biological terms, we define the realized impacts, $\mathcal{I}_T$ and $\mathcal{I}_H$. These metrics integrate the magnitude of the environmental anomaly with the local sensitivity gradient of the persistence threshold
\begin{equation}
    \mathcal{I}_{T}(x)=\left | \frac{\partial \gamma_\text{crit}}{\partial T}\right | |\Delta T(x)|, \quad \mathcal{I}_{H}(x)=\left | \frac{\partial \gamma_\text{crit}}{\partial H^{\max}}\right | |\Delta H(x)|.
\end{equation}

Geometrically, $\mathcal{I}$ represents the first-order displacement of the system toward (or away from) the persistence threshold along each physical axis. This approach separates the magnitude of the climate forcing from the intrinsic response of the local ecosystem.

Numerical mapping of the $\gamma_\text{crit}$ response surface confirms a high structural regularity across the global ocean. The first-order Taylor decomposition remains robust, with a mean $R^2=\GFDLHyperMeanRsq$ and over \GFDLHyperAboveThresholdPerc\% of the global ocean area exhibiting an $R^2>\GFDLHyperThreshold$. This local linearity validates the use of the Thermal Dominance Index ($\mathcal{D}_T$) as a reliable partition of the total climate impact 
\begin{equation}
    \mathcal{D}_T(x)=\frac{\mathcal{I}_T(x)}{\mathcal{I}_T(x) + \mathcal{I}_H(x)}.
\end{equation}

The index $\mathcal{D}_T$ serves as a spatially explicit diagnostic of the regime governing the threshold trajectory
\begin{itemize}
    \item \emph{Thermally Driven ($\mathcal{D}_T>0.5$)}. The trajectory of $\gamma_\text{crit}$ is primarily dictated by thermal-metabolic effects, where temperature-induced cost or gains override physical changes.
    \item \emph{Mixing Driven ($\mathcal{D}_T<0.5$)}. The trajectory of $\gamma_\text{crit}$ is primarily dominated by mixing dynamics, where shifts in MLD -- and its associated effect on light and dilution -- exert the primary control over the persistence threshold.
\end{itemize}

\begin{figure}[tbp]
    \centering
    \includegraphics[width=0.95\textwidth]{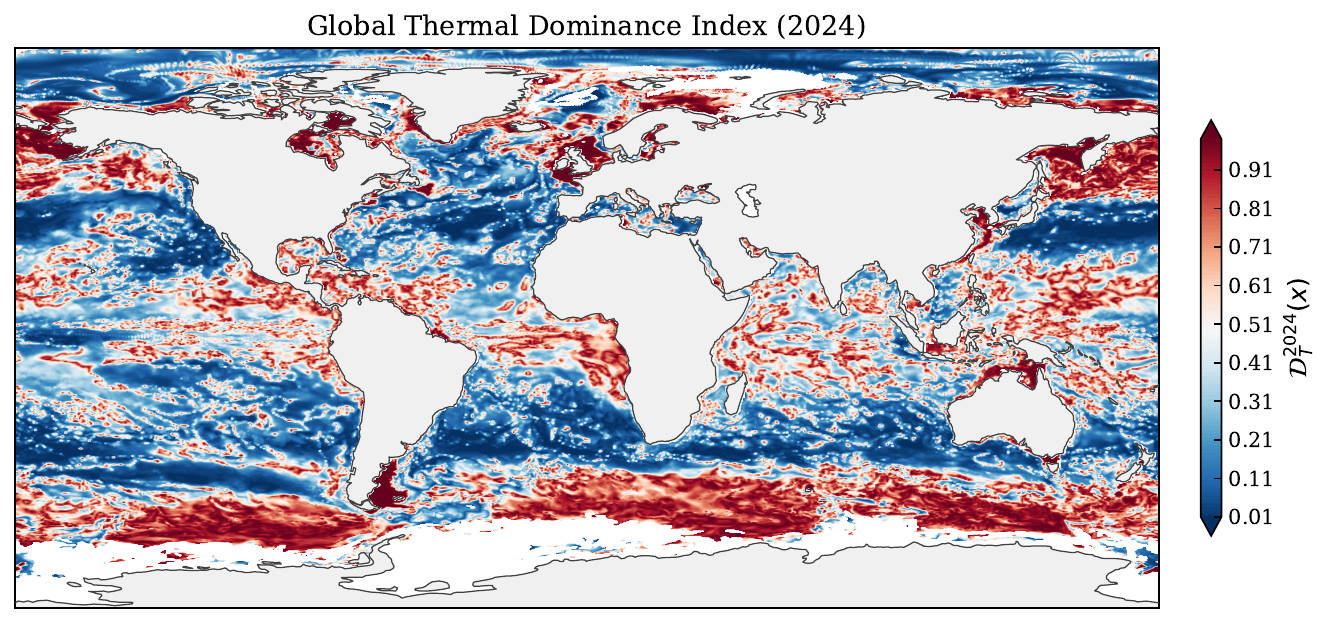}
    \caption{Spatial distribution of the thermal dominance index, $\mathcal{D}_T(x)$. The color scale identifies the primary driver of the persistence threshold: red tones denote thermally driven regimes ($\mathcal{D}_T(x)>0.5$). Blue tones denote mixing driven regimes ($\mathcal{D}_T(x)\le0.5$). The 2024 baseline reveals a heterogeneous ocean where mixing dominance persists in subpolar regions and large oligotrophic gyres.}
    \label{fig:drivers_combined}
\end{figure}

The spatial distribution of $\mathcal{D}_T$ (Figure \ref{fig:drivers_combined}) provides a mechanistic diagnostic for the patterns observed in the persistence threshold (Figure \ref{fig:gamma_crit}). We find that most marginal zones -- where $\gamma_\text{crit}\to 1$ --  are already under strong thermal dominance ($\mathcal{D}_T(x)>0.5$). This is particularly evident in the Southern Ocean, where the realized impact of thermal anomalies consistently overrides the effects of mixing, making the population's persistence highly dependent on the thermal evolution of the water column.

In contrast, the North Atlantic Subpolar Gyre (SPG) represents a different regime of marginality. Although it also exhibits high $\gamma_\text{crit}$ values, the climate impact remains primarily mixing-driven. This suggests that the fate of the SPG is mechanistically constrained by the evolution of its seasonal stratification and the associated dilution stress, rather than by metabolic shifts. This divergence in the dominant control mechanism is fundamental, as it defines which environmental driver is currently pushing these high-latitude ecosystems toward their respective persistence boundaries..

\subsection{Global Projections: Climate-Driven Shifts in $\gamma_\text{crit}$}\label{sec:global_experiments}

The mechanistic drivers identified by the $\mathcal{D}_T$ index lead to significant geographical shifts in the ocean's capacity to support planktonic life. To evaluate how climate change alters global ecosystem viability, we quantify the net displacement of the persistence threshold between the contemporary baseline ($2024$) and the end-of-the-century projections ($2100$). We quantify this shift as
\begin{equation}
    \Delta \gamma(x) = \gamma_\text{crit}^{2100}(x)-\gamma_\text{crit}^{2024}(x)\,.
\end{equation}
For these projections, we employ the SSP5-8.5 high-emission scenario. Our choice is motivated by the objective of characterizing the potential impacts and maximum response of the marine ecosystem under a wide range of anthropogenic forcing, providing a clear benchmark for the physiological pressure exerted by the unmitigated warming and stratification shifts. 

While the continuous shift $\Delta \gamma_{\text{crit}}$ captures the intensity of the physiological pressure, it does not explicitly resolve the crossing of fundamental survival thresholds. To quantify the actual expansion or contraction of habitable niches, we analyze the spatial migration of the critical limit $\gamma_\text{crit}=1$ -- where losses equal gains across the annual cycle -- discretizing the global domain into six distinct ecological regimes

\begin{itemize}
    \item \emph{Stable Viability:} Regions that remain below the extinction threshold ($\gamma_\text{crit} \le1$) in both scenarios, maintaining their capacity for local persistence.
    \item \emph{Stable Restriction:} Regions that remain above the threshold ($\gamma_\text{crit}>1$) in both scenarios.
    \item \emph{Habitat Expansion:} Areas transitioning from restrictive to viable ($\gamma^{2024}_\text{crit}>1 \to \gamma^{2100}_\text{crit}<1$). These represent the colonization of new niches, such as the polar opening, due to the relaxation of physical-metabolic constraints.
    \item \emph{Habitat Contraction:} Areas transitioning from viable to restrictive ($\gamma^{2024}_\text{crit}<1 \to \gamma^{2100}_\text{crit}>1$), where the environment becomes structurally incapable of supporting local populations.
    \item \emph{Ice-Free Viability:} Regions that transition from sea-ice cover to a viable pelagic state ($\gamma^{2024}_\text{crit}=\text{NaN}\to \gamma^{2100}_\text{crit}\le1$).
    \item \emph{Ice-Free Restriction:} Regions that transition from sea-ice cover to a restrictive pelagic state ($\gamma^{2024}_\text{crit}=\text{NaN}\to \gamma^{2100}_\text{crit}\ge1$). While these areas become ice-free, the resulting physical conditions (mixing or temperature) are still too harsh to support a self-sustaining population according to the model.
\end{itemize}

\begin{figure}[tbp]
    \centering
    
    % --- Panel A: Cambio en el Requerimiento Crítico ---
    \begin{subfigure}{\textwidth}
        \centering
        \includegraphics[width=\textwidth]{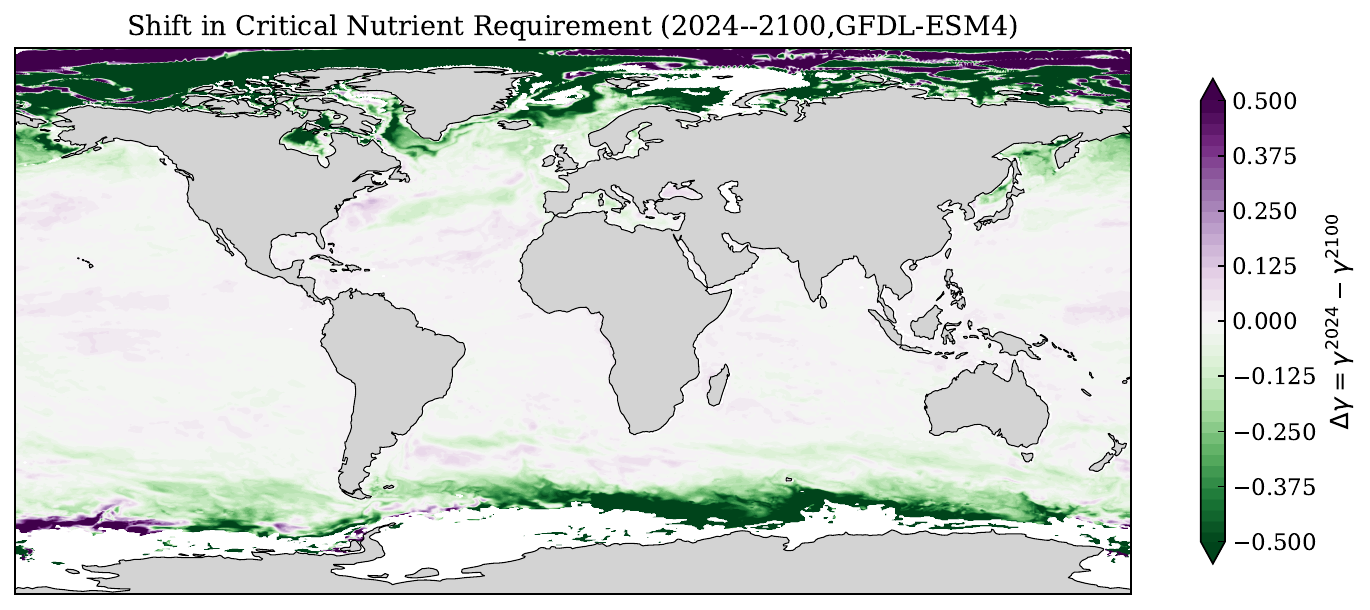}
        \caption{Shift in Critical Nutrient Requirement ($\Delta \gamma = \gamma_{2100} - \gamma_{2024}$)}
        \label{fig:gamma_shift}
    \end{subfigure}
    
    \vspace{0.4cm} 
    
    % --- Panel B: Dominancia Mecanística en el Futuro ---
    \begin{subfigure}{\textwidth}
        \centering
        \includegraphics[width=\textwidth]{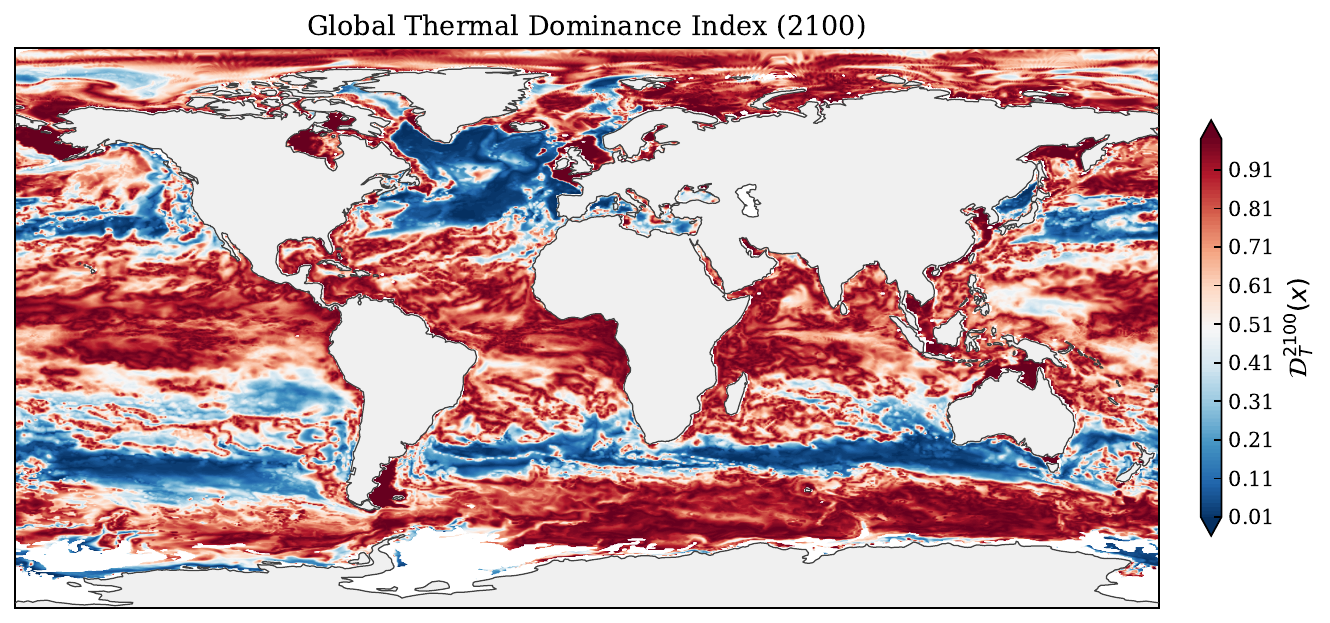} 
        \caption{Environmental Control Projection (2100)}
        \label{fig:dominance_2100}
    \end{subfigure}
    
    \caption{(a) Critical Nutrient Shift ($\Delta \gamma_\text{crit}$): Changes in the persistence threshold between 2024 and 2100 under the high-emission SSP5-8.5 scenario. Green tones ($\Delta \gamma < 0$) indicate lower critical nutrient requirements, while purple tones ($\Delta \gamma > 0$) indicate greater nutrient requirements. 
    (b) Thermal Dominance Index: Spatial distribution of the thermal dominance index $\mathcal{D}_T$ in $2100$. The widespread expansion of red tones ($\mathcal{D}_T > 0.5$) reflects the ``metabolic tropicalization'' of high latitudes.}
    \label{fig:mechanisms_shift}
\end{figure}

\begin{figure}[tbp]
    \centering
    \includegraphics[width=\textwidth]{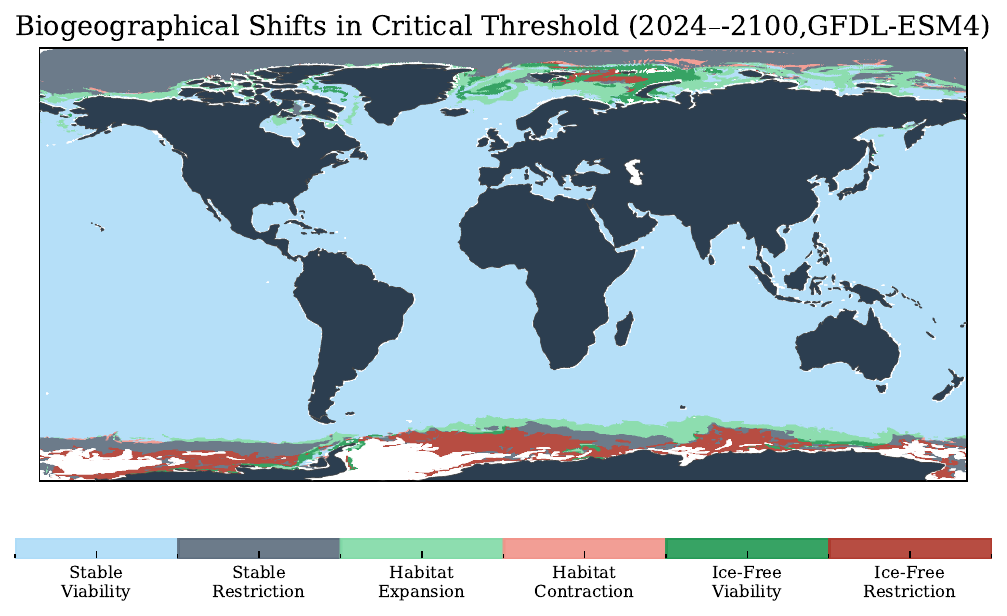}
    
    \caption{Spatial categorization of ecological transitions based on the $\gamma_\text{crit}=1$ threshold. 
    Stable regimes are shown in ultra-light blue (Viable) and blue-grey (Restrictive). 
    Transitions in historically open waters include Habitat Expansion (vivid green) and Contraction (alert red). 
    Cryospheric retreat reveals an asymmetric Polar Opening: New Ice-Free Habitats (intense green) represent successful niche expansions, primarily in the Arctic, while New Ice-Free Restrictive areas (dark red) indicate regions where geographical opening is countered by extreme physical constraints, predominantly in the Southern Ocean.}
    \label{fig:habitat_transitions}
\end{figure}

The global analysis (Figure \ref{fig:gamma_shift}) reveals a dominant trend of critical nutrient requirement reduction ($\Delta \gamma_\text{crit} <0$) across most high-latitude regions, indicating that the nutrient threshold for persistence becomes less demanding by $2100$. This biological response is mechanistically underpinned by the global shift in physical forcing (Figure \ref{fig:dominance_2100}), where a widespread increase in thermal dominance ($\mathcal{D}^{2100}_T>0.5$) indicates a transition towards metabolically-controlled regimes. Under this high-forcing scenario, small regions of low increased nutrient requirements emerge ($\Delta \gamma_\text{crit} >0$) at low-mid latitudes where temperature also dominates, indicating its dual role: favoring the acceleration of metabolic rates and also metabolic stress.

This trend is not uniform. At the extreme polar margins localized bands of critical nutrient requirement increase ($\Delta \gamma_\text{crit} >0$) emerge. The driver analysis (Figure \ref{fig:drivers_combined} and \ref{fig:dominance_2100}) reveals that these regions transition toward a thermally-dominated regime by $2100$, suggesting that thermal limitation remains the primary constraint. Even as sea-ice retreats, the exposed waters can remain at temperatures low enough to suppress metabolic rates, increasing the nutrient requirements for persistence. These ``cold-traps'' -- along with small spots where extreme warming might trigger thermal stress -- counter the general reduction of the high latitudes.

\begin{table}[htbp]
    \centering
    \caption{Global area and ocean coverage of phytoplankton persistence regimes under the GFDL-ESM4 SSP5-8.5 projection (2015--2100).}
    \label{tab:habitat_stats_gfdl}
    \begin{tabular}{l c c}
    \toprule
    \textbf{Ecological Regime / Habitat} & \textbf{Area ($10^6$ km$^2$)} & \textbf{Ocean Coverage (\%)} \\
    \midrule
    Stable Viability & \GFDLStableViabilityArea & \GFDLStableViabilityPerc \\
    Stable Restriction & \GFDLStableRestrictionArea & \GFDLStableRestrictionPerc \\
    Habitat Expansion & \GFDLHabitatExpansionArea & \GFDLHabitatExpansionPerc \\
    Habitat Contraction & \GFDLHabitatContractionArea & \GFDLHabitatContractionPerc \\
    Ice-Free Viability & \GFDLIceFreeViabilityArea & \GFDLIceFreeViabilityPerc \\
    Ice-Free Restriction & \GFDLIceFreeRestrictionArea & \GFDLIceFreeRestrictionPerc \\
    \bottomrule
    \end{tabular}
\end{table}

To determine the ecological impact of these shifts, we evaluate the resulting ecological transitions (Figure \ref{fig:habitat_transitions}, Table \ref{tab:habitat_stats_gfdl}). Quantitatively, the results reveal a global trend toward niche expansion: while Habitat Contraction is almost negligible ($\GFDLHabitatContractionArea$  \unit{million.km^2}), Habitat Expansion into previously restrictive open waters covers $\GFDLHabitatExpansionArea$ \unit{million.km^2}. This expansion appears as distinct bands in both hemispheres, where the relaxation of metabolic-physical constraints promotes viability; a ``Polar Opening''.

However, a fundamental disparity emerges in the nature of the ``Polar Opening''. Although cryospheric retreat reveals vast new territories, only $\GFDLIceFreeViabilityArea$ \unit{million.km^2} are Ice-Free Viability. In contrast, a much larger area -- approximately $\GFDLIceFreeRestrictionArea$ \unit{million.km^2} -- remains as Ice-Free Restriction. These figures confirm that geographical ice retreat is a necessary but not sufficient condition for niche expansion. For every $1$ \unit{km^2} of viable habitat gained by the retreat of the ice, nearly $4$ \unit{km^2} remain structurally hostile, particularly in the Southern Ocean. This underscores that while the anthropogenic forcing generally promotes critical nutrient requirement reduction in historically open waters, the newly exposed polar frontiers remain largely locked under extreme environmental constraints.

\subsubsection{The North Atlantic Subpolar Gyre (SPG): A Metabolic Refuge}

Beyond the broad latitudinal trends, the global maps reveal the SPG as a unique regional exception. This area is characterized by the ``warming hole'' phenomenon -- a persistent anomaly where surface temperatures remain nearly constant or even decrease -- deviating from the global warming trend \cite{rahmstorf_is_2024}. Under the high-emission forcing of the SSP5-8.5 scenario, this region introduces an important tension between the mixing dynamics and the anomalous temperature profile that we can study based on our stability metric $\gamma_\text{crit}$ and our diagnostic profile $\mathcal{D}_T$.
\begin{figure}[tbp]
    \centering
    \begin{subfigure}{\textwidth}
        \centering
        \includegraphics[width=\textwidth]{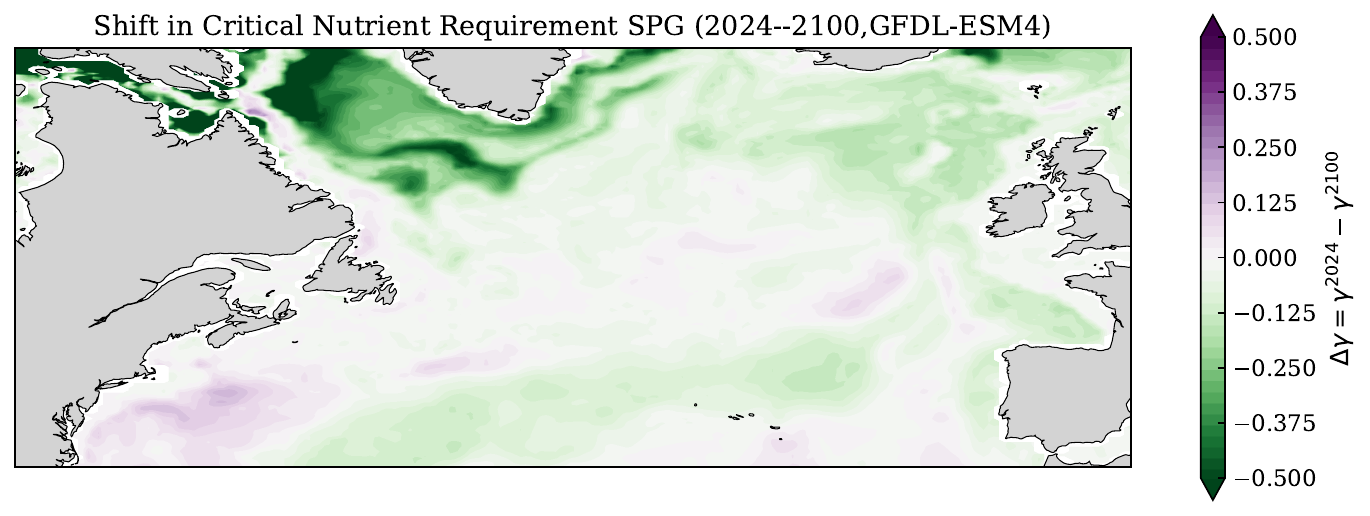} 
        \caption{Regional Shift in Critical Nutrient Requirement ($\Delta \gamma$)}
        \label{fig:spg_nutrient_shift}
    \end{subfigure}
    
    \vspace{0.4cm} 
    
    \begin{subfigure}{\textwidth}
        \centering
        \includegraphics[width=\textwidth]{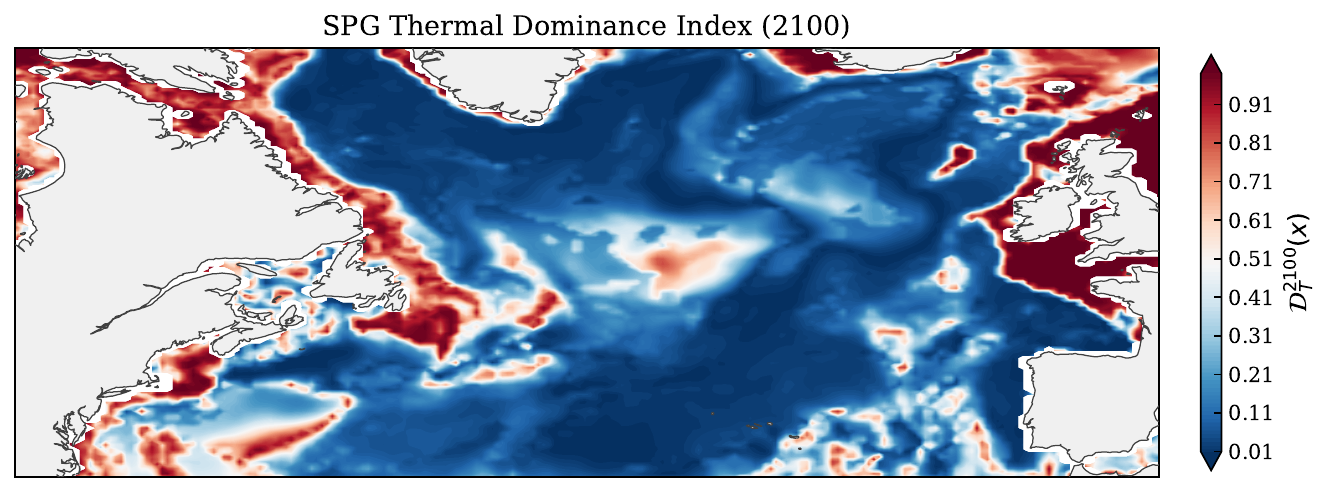} 
        \caption{Environmental Control Projection (2100)}
        \label{fig:spg_control_2100}
    \end{subfigure}

 \caption{ (a) Regional Nutrient Shift: The core of the SPG exhibits a dampened response (pale tones), contrasting with the intense reduction in nutrient requirements (deep green) observed in the surrounding boundary currents. 
    (b) Regional Thermal Dominance Index (2100): Spatial distribution of the thermal dominance index $\mathcal{D}_T$. While the North Atlantic basin largely transitions toward a thermally-controlled regime ($\mathcal{D}_T > 0.5$, red), the SPG remains a persistent mixing-driven domain (blue tones, $\mathcal{D}_T < 0.3$).}
    \label{fig:spg_analysis}
\end{figure}

As shown in Figure \ref{fig:spg_nutrient_shift} the core of the SPG exhibits a dampened biological response. While the surrounding boundary currents undergo a rapid reduction in nutrient requirements (deep blue, $\Delta \gamma < 0$), the Gyre's core shows only a faint shift ($\Delta \gamma \approx  0$). This relative stability is remarkable given the high-emission forcing of the SSP5-8.5 scenario, and it is mechanistically underpinned by the persistence of the mixing-driven regime. Our driver control analysis (Figure \ref{fig:spg_control_2100}) confirms that, while the high latitudes of the North Atlantic basin shifts toward thermal dominance (red tones, $\mathcal{D}_T>0.5$), the SPG core remains one the few regions where mixing constraints continue to govern the persistence threshold by 2100 (blue domain, $\mathcal{D}_T<0.3$ ). 

This suggests that the SPG may act as a metabolic refuge, since regional thermal forcing remains relatively stable, the ecosystem is partially shielded from the ``metabolic relief'' observed in warmer basins. Rather than a transition toward a new state, the SPG represents a case of environmental anchoring, where historical physical constraints -- both thermal and mixing -- prevail against the encroachment of ``global tropicalization.''

\section{Discussion}

We structure the discussion to first address the computational and methodological advantages of our framework. Subsequently, we examine the ecological and physical implications of our results, focusing on three key advances: (1) a quantitative re-evaluation of the hierarchy of physical drivers governing phytoplankton dynamics, (2) the identification of global trends, such as the global ecological expansion, the asymmetric polar opening and the SPG as a metabolic refuge, and (3) the mathematical implications of our stability theorems for the persistence of marine ecosystems.

\subsection{A General Approach for Sensitivity of Ecosystem Dynamics to Model Assumptions}

Earth System Models (ESMs) face a fundamental trade-off between biological realism and computational feasibility. Increasing ecosystem complexity not only demands greater processing power but also introduces structural uncertainties that obscure the causal link between physical forcing and biological rates \cite{anderson_2005}, a challenge that persists even in the latest CMIP6 generation where increased complexity has not necessarily reduced projection uncertainty \cite{kwiatkowski, heneghan_2021}.

Here we demonstrate an approach to test the impacts of modifications to the parameterization of ecosystems in a climate model that requires minimal computational expense. By leveraging the linearization around the extinction equilibrium ($E_0$), we can analytically derive stability metrics such as the critical nutrient requirement ($\gamma_\text{crit}$) and thermal dominance index ($\mathcal{D}_T$) directly from physical outputs. This approach provides a reproducible tool for rapidly evaluating impacts of climate changes on ecosystems in a way that incorporates novel understanding of ecosystem processes \cite{martin_when_2024}. Whereas standard outputs (e.g., chlorophyll) typically show the net result of competing stressors, our metrics explicitly disentangle the mechanistic contributions of thermal-metabolic forcing versus mixing constraints. This transparency is important for identifying regional anomalies, such as the North Atlantic ``warming hole'', where the physical drivers diverge from global trends.

\subsection{Beyond the Critical Depth: The Emerging Dominance of Metabolic Constraints}
Our global analysis of the thermal dominance index ($\mathcal{D}_T$) provides a new perspective on the universality of the classical CDH \cite{sverdrup}. The CDH has previously faced scrutiny regarding the role of predator-prey decoupling \cite{behrenfield2010} and the physical mechanism of turbulence shutdown \cite{taylor_ferrari_2011}. Our approach unifies these critiques, which have previously been studied in terms of bloom timing, with a fundamental \emph{metabolic} extension to this framework. We introduce a metric that generalizes the critical depth approach to multiple physical stressors and ecosystem processes. We suggest that the classical mixing-light balance is intrinsically linked to a temperature-dependent metabolic demand that scales with the thermal state of the environment.

By utilizing the SSP5-8.5 scenario as a systemic stress test, we observe that the regime described by the CDH may be increasingly limited to the North Atlantic SPG in future climate scenarios. The 2100 projections (Figure \ref{fig:dominance_2100}) reveal a significant geographical expansion of metabolic control ($\mathcal{D}_T>0.5$), where the acceleration (or deceleration) of biological rates becomes the primary driver of the persistence threshold. In this context, the metabolic sensitivity to thermal shifts acts as a dominant constraint that can either override the traditional advantages of a shallower MLD or exacerbate the stress of deep mixing.

We propose that the $\mathcal{D}_T$ index serves as a geographical validity map for the CDH. In regions where mixing remains the primary driver ($\mathcal{D}_T\approx0$), Sverdrup's original formulation continues to be a robust predictive tool. However, as the ocean undergoes thermal restructuring, these areas are projected to contract, suggesting that the stability of marine ecosystems will be increasingly governed by the dynamic competition between environmental energy supply and the internal metabolic requirements.

\subsection{Ecological Expansion and the Polar Opening Paradox}

The global reorganization of marine ecosystems under climate change is frequently characterized as a poleward migration of productive niches \cite{boyce}. Our analysis of the shift in the critical nutrient requirement ($\Delta \gamma_\text{crit}$) and the resulting ecological transitions confirms this trend, but reveals a complex balance between metabolic relief, habitat expansion, and physical structural barriers.

First, our results could imply a thermodynamic basis for the declines in productivity observed in low-to-mid latitudes \cite{hong}. We identify critical zones where temperature leads to an increase in the critical nutrient requirement ($\Delta \gamma_\text{crit}>0$). In these regions, temperature acts as a metabolic burden: the acceleration of biological rates raises the ``metabolic tax,'' forcing the population to require higher nutrient saturation levels merely to maintain the extinction-persistence equilibrium. This suggests that, in specific low-latitude regimes where populations are not light-limited, the metabolic cost of warming occurs without benefits of increased stratification, pushing the ecosystem closer to extinction thresholds.

Second, our results quantify the extent of the high-latitude ``greening'' \cite{polovina_2011, ardyna_arrigo_2020}. Under the sustained forcing of the SSP5-8.5 scenario, we identify a massive Habitat Expansion of approximately $8.6$ \unit{million.km^2}, where the critical nutrient threshold decreases. This expansion, organized in two coherent circumpolar bands, suggests that ``metabolic tropicalization'' fundamentally lowers the barrier for phytoplankton persistence across vast new areas. 

However, this expansion is subject to a physical-biological decoupling at the poles, which our stress-test approach reveals as a ``Polar Opening Paradox''. While the net balance of tropicalization is positive, the efficiency of the ice retreat is low. We find that for every square kilometer where the threshold becomes viable due to ice retreat, nearly four square kilometers remain locked in Ice-Free Restriction ($7.58$ \unit{million.km^2}). This $1:4$ ratio, however, is not symmetrically distributed between the hemispheres, reflecting a profound polar asymmetry in the ``Opening''. While the Arctic experiences a widespread ``greening'' as nutrient requirements decrease and sea-ice retreat aligns with metabolic relief, the Southern Ocean remains a bastion of Ice-Free Restriction. In the Southern Hemisphere, the rigid mixing-light regime -- even under extreme warming -- acts as a physical ceiling. This suggests that while the North Pole may transition to a more productive state, ice retreat in the Antarctic is a necessary but not sufficient condition for niche expansion. Thus, the newly exposed waters may remain low biomass in current models. However, observations of under-ice blooms even in the Southern Ocean suggest that modeling the unique metabolic capabilities of high-latitude communities is likely to be informative \cite{horvat_evidence_2022,mcclish_majority_2023}.

\subsection{The North Atlantic Subpolar Gyre (SPG): A Metabolic Refuge in a Changing Ocean}

An exception to the widespread global decline in the critical nutrient requirement is the SPG, specifically the region associated with the North Atlantic ``warming hole'' or cold anomaly. While this feature is widely attributed to the slowdown of the Atlantic Meridional Overturning Circulation \cite{rahmstorf_2015, keil_2020}, its biological implications have largely been viewed through the lens of nutrient transport or physical stratification.

Our stability analysis (Figure \ref{fig:spg_nutrient_shift}) offers a novel perspective: the ``warming hole'' acts as a \emph{metabolic refuge}. Unlike the majority of the global ocean, where the ``metabolic tax'' of temperature forces a reorganization of persistence thresholds, the SPG ecosystem remains anchored in a low-cost thermodynamic regime. Our results show that in this region, the critical nutrient requirement exhibits a relative stasis ($\Delta \gamma_\text{crit}\approx 0$). Even under the high-emission SSP5-8.5 scenario, this stability indicates that the SPG is shielded from the global trend of metabolic acceleration, allowing the ecosystem to maintain its historical balance between nutrient supply and biological demand.

As shown in our driver dominance analysis (Figure \ref{fig:dominance_2100}), the SPG remains a mixing-driven bastion ($\mathcal{D}_T<0.3$) in a world that is otherwise transitioning toward thermal control. This ``environmental anchoring'' means that the local phytoplankton populations do not experience the same pressure to ``run faster to stay in place'' that characterizes the tropical and sub-tropical oceans.

The existence of this metabolic refuge emphasizes the importance of regional circulation dynamics in modulating global metabolic trends. If the physical processes that maintain the cold anomaly persist -- even as anthropogenic forcing intensifies -- the SPG could serve as a vital reservoir of biological stability. It preserves a metabolic signature that is rapidly disappearing elsewhere in the 21st-century ocean, suggesting that regional oceanographic ``anomalies'' may become critical buffer zones for marine biodiversity.

\subsection{From Local Instability to Global Persistence}

The classical understanding of phytoplankton dynamics has long been centered on the conditions for bloom initiation. Sverdrup's Critical Depth Hypothesis \cite{sverdrup} fundamentally treats the bloom as a local instability, where growth temporarily escapes loss. However, this perspective has faced significant criticism for its inability to account for the seasonal decoupling of grazing and division rates, as well as the maintenance of populations during the ``non-blooms'' periods \cite{behrenfield2010}. Our results resolve this tension by shifting the focus from transient instabilities to the global attractor of the system and the conditions for uniform persistence \cite{zhao, wang2008threshold}.

By treating the marine ecosystem as a non-autonomous dynamical system, we move beyond the instantaneous balances of traditional models. Our stability criterion -- formally the Floquet exponent $\lambda_P$-- quantifies the invasion growth rate, identifying whether the extinction equilibrium $E_0$ is a repeller for the seasonal flow. This aligns our methodology with Modern Coexistence Theory \cite{Chesson2000}, and the framework of periodic persistence \cite{zhao}, where the fundamental metric for population maintenance in fluctuating environments is the ability to recover from seasonal minima. The critical threshold $\gamma_\text{crit}$ is not just a seasonal trigger for a bloom, it is a formal bifurcation point that separates a regime of permanent extinction from one of global persistence.

This mathematical bridge allows us to interpret 21st-century climatic shifts with a new level of rigor. The widespread ecological transitions identified in high-latitudes (see Figure \ref{fig:habitat_transitions}) represent a geographic displacement of these bifurcation points. Under the stress-test of the SSP5-8.5 scenario, we observe that sea-ice retreat and metabolic shifts do more than alter the timing of seasonal events; they redefine the boundaries of the niche itself. By framing climate change as a global shift in bifurcation points, we provide a formal basis for understanding how the ``metabolic tax'' and physical constraints determine the functional integrity of the marine biosphere.

\subsection{Limitations and Future Directions}

While our framework offers a robust diagnostic for marine ecosystem stability, the assumptions inherent in our modeling choices suggest critical avenues for future research.

First, we distinguish between biological traits and biological rates. In this study, the physiological traits defining the TPC -- such as the optimal temperatures -- were held constant to represent a generic plankton functional type. Thus, our results isolate the thermodynamic pressure exerted by the environment from the compensatory mechanism of evolutionary adaptation. Future work should relax the assumption of fixed traits by incorporating adaptive dynamics or evolutionary game theory \cite{Litchman2007,walworth_microbial_2020}. This would allow for an exploration of how survival strategies and trait distributions shift in response to climate change, moving from a question of mere persistence to one of community composition and niche construction.

Second, while we incorporated the temperature dependence of metabolic rates via the EEP-formulation \eqref{eqn:exp-eep}, our analysis focuses exclusively on the persistence of the primary producer. However, evidence suggest that the predator-prey parameters, such as handling time and search efficiency within the Holling functional response, are also temperature-dependent \cite{ZHAO2020}. 

Integrating these thermally-sensitive functional responses is a necessary next step to examine the stability and persistence of higher trophic levels. By applying our non-autonomous framework to a full NPZ-T interaction, we could diagnose regions where climate change triggers a metabolic mismatch between producers and consumers. Such a mismatch can lead to the emergence of unstable, high-amplitude oscillations or other disruptive dynamics.

We expect that the results are sensitive to the structure of the ecological model used. This reduced modeling approach should be applied to other models and may be useful for diagnosing the mechanistic drivers of ecosystem change. Extending our $\gamma_\text{crit}$ analysis to include more complex ecosystems with differential responses to temperature multi-trophic feedbacks would allow us to map where the projected ``metabolic tax'' of a changing ocean may compromise the seasonal coupling of the food web, transitioning the ecosystem from balanced annual cycle toward a regime of unpredictable ecological disruptions. It could be particularly relevant to study complex ecosystems with strain level metabolic diversity in this framework \cite{krinos_intraspecific_2025}.

\section{Conclusion}

This study proposes a theoretical framework to evaluate phytoplankton stability that transcends instantaneous balances, focusing on the conditions for long-term persistence throughout the annual cycle. By applying non-autonomous dynamical systems theory to 21st-century climate projections, we generalize previous theories and incorporate thermodynamic factors.

Our results, characterized under a high-emission scenario, lead to three main conclusions. First, we observe a global trend toward a reduction in the critical nutrient requirement ($\Delta \gamma_\text{crit}<0$), driven by a widespread transition toward metabolic control. Second, we reveal a significant decoupling between physics and biology at the poles; the 1-to-4 ratio between viable niches and ice-free deserts demonstrates that the retreat of the cryosphere does not guarantee a proportional expansion of life. Finally, we identify the North Atlantic Subpolar Gyre as a metabolic refuge where mixing dynamics anchor the ecosystem, shielding it from global metabolic shift, maintaining stability conditions that are disappearing elsewhere.

Beyond these regional shifts, this framework provides a projection of the future ocean's complexity. Through the thermal dominance index ($\mathcal{D}_T$), we are able to deconstruct whether an ecosystem's fate is driven by physical mixing or metabolic costs. Ultimately, our methodology offers a mechanistic basis to understand that the ocean's response to extreme anthropogenic forcing is not uniform, but the result of a dynamic balance between environmental energy supply and the escalating metabolic demands of life.

\appendix

\section{Model Well-posedness}\label{sec:model-well-posedness}

\counterwithin{equation}{section}

\begin{proposition}\label{th:well_pos}
Given any initial condition $(N(0), P(0), Z(0)) \in \mathbb{R}^3_+$, the system \eqref{eqn:NPZ-T-full} admits a unique global solution. The non-negative orthant $\mathbb{R}^3_+$ is positively invariant. Moreover, the total mass $C(t) = N(t) + P(t) + Z(t)$ is conserved for all $t \geq 0$, and the solutions are confined to the compact simplex:
\begin{equation*}
X = \{ (N,P,Z) \in \mathbb{R}^3_+ : N + P + Z = C_0 \}\,.
\end{equation*}
\end{proposition}

\section{Mathematical Foundations: Floquet Theory and Uniform Persistence}\label{sec:foundations}

To analyze the stability of the NPZ-T model, we treat it as a non-autonomous periodic system. Since the coefficients (temperature and mixing) satisfy $A(t)=A(t+T_p)$ with $T_p=365\unit{days}$,  classical autonomous stability is insufficient.

\subsection{Local Stability and the Monodromy Matrix}

The stability of the extinction equilibrium $E_0 = (C_0,0,0)$ -- where both Phytoplankton ($P$) and Zooplankton ($Z$) are absent -- is determined by the linearization of the system around the origin. We define the monodromy matrix $C$ as the operator that maps the initial perturbation $\Phi(0)=(P(0),Z(0))$ to its state after one full cycle
\begin{equation*}
\Phi(T_p) = C \Phi(0)\,.
\end{equation*}

The eigenvalues of $C$, known as Floquet multipliers ($\rho$), allow us to define the spectral radius $r(C) = \max |\rho|$, that dictates the local fate of the system
\begin{itemize}
\item If $r(C) < 1$, $E_0$ is locally asymptotically stable, \emph{i.\,e.}, any small population perturbation in a neighborhood $U \subset \mathbb{R}^2_+$ centered at $E_0$ will decay to zero ($\lim_{t \to \infty} \Phi(t) = 0$).
\item If $r(C) > 1$, $E_0$ is locally asymptotically unstable.
\end{itemize}

\subsection{From Instability to Uniform Persistence}

While local instability suggests growth near the origin, the long-term survival of the ecosystem requires uniform persistence. This means there exists an $\epsilon > 0$ such that for any initial condition with $P(0), Z(0) > 0$
\begin{equation*}
\liminf_{t \to \infty} P(t) \geq \epsilon, \quad \liminf_{t \to \infty} Z(t) \geq \epsilon.
\end{equation*}

Our analytical framework (Theorem \ref{th:persistence}) leverages the fact that for this periodic NPZ structure, the instability of $E_0$ (characterized by a positive Floquet exponent $\lambda_P = \frac{1}{T_p} \ln r(C) > 0$) is sufficient to guarantee uniform persistence for the phytoplankton population. In this regime, the extinction state acts as a repeller, ensuring the ecosystem remains bounded away from the boundary of the state space.

\subsection{Key References}

\begin{itemize}
\item \emph{Floquet Theory:} \cite{Hale} for linear periodic systems and monodromy matrices and \cite{klausmeier_floquet_2008} for useful examples applied to ecology.
\item \emph{Persistence in Periodic Models:} \cite{zhao, wang2008threshold} for the formal link between the spectral radius of the periodic operator and global ecosystem persistence.
\end{itemize}

\section{Proofs of Theorems}\label{sec:proofs}

\subsection{Proof of Proposition \ref{th:well_pos}}

To establish the well-posedness of the model, we verify the existence, uniqueness, positivity, and conservation of mass of the solutions. Let $\Phi = (N, P, Z)^\top \in \mathbb{R}^3$ be the state vector. The system \eqref{eqn:NPZ-T-full} can be written as the Cauchy problem
\begin{equation} \label{eqn:cauchy}
    \frac{d \Phi}{d t} = F(t, \Phi), \quad \Phi(0) = \Phi_0 \in \mathbb{R}^3_+\,.
\end{equation}
The functions involved in the vector field $F(t, \Phi)$ are $C^1$ and bounded on the domain of interest. Therefore, $F$ is locally Lipschitz continuous, and by the Picard-Lindelöf Theorem, there exists a unique local solution $\Phi(t)$ on $[0, T_{\max})$.

Regarding positivity, the planes $P=0$ and $Z=0$ are invariant sub-spaces of the system, since $P=0 \implies \dot{P}=0$ and $Z=0 \implies \dot{Z}=0$ (noting that $h(t,0,T)=0$). Also, on the boundary $N=0$, we have:
\begin{equation*}
\left. \frac{d N}{d t} \right|_{N=0} = (1-\alpha)h(t,P,T) Z + m_P(t,T) P + m_Z(t,T) Z^2 + s_+(t)(P+Z) \ge 0\,,
\end{equation*}
for all $(P,Z) \in \mathbb{R}^2_+$. Thus, the non-negative orthant $\mathbb{R}^3_+$ is positively invariant.

Adding the equations, we find $\frac{d}{dt}(N+P+Z) = 0$, which implies $C(t) = C_0$ for all $t \ge 0$. This conservation law, combined with the positivity of the variables, ensures that the solutions are confined to the compact simplex:
\begin{equation*}
X = \{ (N,P,Z) \in \mathbb{R}^3_+ : N+P+Z = C_0 \}\,.
\end{equation*}
Hence, the solutions are bounded, and no finite-time blow-up occurs, proving that $T_{\max} = \infty$. \qed

\subsection{Proof of Theorem \ref{th:extinction}}

We establish global extinction conditions using standard comparison arguments. First, we construct a differential inequality for the phytoplankton equation $P(t)$. From Proposition \ref{th:well_pos}, we know that the total mass is conserved, so $N(t) \le C_0$ for all $t\ge0$, then the nutrient limitation term satisfies $g(N) \le g(C_0)$. So, we obtain
\begin{align*}
    \frac{d P}{dt} & = f(t,I,H)g(t,N,T) P - h(t,P,T) Z - m_P(t,T) P - s_+(t) P \\
                     & \le \left( f(t,I,H)g(t,C_0,T) - m_P(t,T) - s_+(t) \right) P\,.
\end{align*}
Let $P_l(t)$ be the solution to the following linear comparison system
\begin{equation*}
    \frac{d P_l}{dt} = \left( f(t,I,H)g(t,C_0,T) - m_P(t,T) - s_+(t) \right) P_l\,.
\end{equation*}
A standard comparison argument shows that as $P(0) = P_l(0)$, then $0 \le P(t) \le P_l(t)$ for all $t\ge0$.
Since the coefficients are $T_p\text{-periodic}$ we can write 
\begin{equation*}
    P_l(nT_p)=P(0)e^{nT_p \lambda_P}.
\end{equation*}
As $n\to \infty$, $Pl(nT_p)\to 0$, since $ \lambda_P < 0$. The per-capita growth rate of $P$ is bounded from above, letting it be $M$
\begin{equation*}
\frac{\dot{P}}{P}\leq M\,.
\end{equation*}
Let $n\in \mathbb{N}$, the above implies that for all $t\in\left[nT_p,nT_p+T_p\right]$ we get
\begin{equation*}
P_l(t)\leq P(nT_p)\exp(M(t-nT_p))\leq P(nT_p)\exp(MT_p)\,.
\end{equation*}
As $n\to \infty$, implying that $\lim_{t \to \infty} P_l(t) = 0$. Thus
\begin{equation*}
    \lim_{t \to \infty} P(t) = 0\,.
\end{equation*}

Recall the dynamics of $Z$ 
\begin{equation*}
    \frac{d Z}{dt} = \alpha h(t,P,T) Z - m_Z(t,T) Z^2 - s_+(t) Z\,.
\end{equation*}
Since $h(t,P,T)$ is continuous and $h(t,0,T)=0$, the fact that $P(t)\to 0$ as $t\to \infty$ implies that for any $\eta > 0$, there exists a time $t_0$ such that $\alpha h(t,P(t), T) < \eta$ for all $t > t_0$. Choosing $\eta < \tfrac{1}{T_p}\int_{t_0}^{t_0+T_p}s_+(\tau)d\tau$, we obtain the inequality
\begin{equation*}
    \frac{d Z}{dt} \le \left( \eta - s_+(t) \right) Z -m_Z(t,T) Z ^2 \le \left( \eta - s_+(t) \right) Z\,,
\end{equation*}
Integrating this inequality, we have
\begin{equation*}
Z(t) \le Z(t_0) \exp \left( \int_{t_0}^t (\eta - s_+(\tau)) d\tau \right)=Z(t_0)F(t),~~t\ge t_0\,.
\end{equation*}
As $F(t_0+T_p)<1$, because of the choose of $\eta$, and by a similar argument as in the case of phytoplankton, this implies that 
\begin{equation*}
    \lim_{t \to \infty} Z(t) = 0\,.
\end{equation*}

Since $P(t) \to 0$ and $Z(t) \to 0$ as $t\to \infty$, by conservation of the mass (Proposition \ref{th:well_pos}), the nutrient concentration must approach the total inventory $N(t) \to C_0$. Thus, the system globally converges to the extinction equilibrium $E_0 = (C_0, 0, 0)$.
\qed

\subsection{Proof of Theorem \ref{th:persistence}}

To prove this theorem, we will use Poincaré map theory applied to non-autonomous periodic dynamical systems. We first recall the state space of the system $X:=\{(N,P,Z)\in \mathbb{R}^3:N+P+Z=C_0\}$ given by Proposition \ref{th:well_pos}, that is a compact simplex. We then define the Poincaré map $f:X\to X$ given by $f(x_0)=\Phi(T_p,x_0)$, where $\Phi$ comes from the Cauchy problem \ref{eqn:cauchy}. The map $f$ is continuous on $X$ due to the continuous dependence of the solutions with respect to the initial conditions. Also, we define $X_0:=\{(N,P,Z)\in X: P>0\}$ as the phytoplankton persistence space.

We will first show that the Poincaré map $f$ is uniformly persistent using Theorem 1.3.3 from Zhao \cite{zhao}. First, we note that $f(X_0)\subset X_0$. This is obtained directly from the fact that the phytoplankton equation can be written as $\dot{P}=P \mathcal{H} (\Phi(t,x_0))$, where $\mathcal{H}$ is a continuous function. By the analytical expression of the solution $P$, any trajectory with $P(0)>0$ remains strictly positive for any finite time, so $f^n(X_0)\subset X_0~~\forall n\in \mathbb{N}$. Then, given that the state space $X$ is compact and forward-invariant under $f$, it follows that $f$ is compact and point dissipative. Thus, by Theorem 1.1.3 from \cite{zhao}, there exists a global attractor $A\subset X$. 

Let us see that $f$ is weakly uniformly persistent with respect to $(X_0,\partial X_0)$ via contradiction. Let's assume that $\exists x\in X_0$ such that 
\begin{equation*}
\lim_{n\to \infty} P_{x_0}(nT_p)=0\,.
\end{equation*}
Analogous to what was demonstrated in the Theorem \ref{th:extinction}, it follows that the convergence of the discrete sequence $P_{x_0}(nT_p)\to 0$ implies $P_{x_0}(t)\to 0$ as $t\to \infty$, therefore, $E_0$ is globally asymptotically stable. Thus, the $\omega\text{-limit}$ set of any such trajectory is the singleton $\{E_0\}$.

Now let's study the behavior near the extinction equilibrium. First, we note that the Poincaré map of the $P$ component is given by
\begin{equation*}
P_{n+1}=P_n \exp\left( \int_0^{T_p} \mathcal{H}(\Phi(t,x_n)) dt\right):=P_n G(x_n)\,,
\end{equation*}
where
\begin{equation*}
\mathcal{H}(t,N,P,Z)=f(t,I,H)g(t,N,T)-\tfrac{h(t,P,T)Z}{P}-m_P(t,T)-s_+(t)
\end{equation*}
is the per-capita growth rate of the phytoplankton. The continuity of $\mathcal{H}$ at $E_0$ follows from the fact that the grazing term $h(t,P,T)Z$ is of class $\mathcal{C}^1$ with $h(t,0,T)=0$, ensuring that the per-capita loss $h/P$ vanishes as $P \to 0^+$.

Since $\lambda_P>0$, it follows that $G(E_0)=\exp(T_p\lambda_P)>1$. Let $\sigma$ be a constant such that $1<\sigma<G(E_0)$, and we define $\epsilon=G(E_0)-\sigma>0$. By the continuity of $\mathcal{H}$ and the flow, there exists a neighborhood $\mathbb{B}_{\delta}(E_0)$ with $\delta>0$ such that $G(x)>\sigma$ for all $x\in \mathbb{B}_{\delta}(E_0)$. By convergence to $E_0$, there exists some $N_0$ such that the trajectory $(P_{N_0+n})_n$ remains inside the ball $\mathbb{B}_{\delta}(E_0)$ for all $n \ge 0$, therefore
\[
P_{N_0+k}=P_{N_0+k-1}G(x_{N_0+k-1})\geq P_{N_0+k-1}\sigma\geq \dots \geq P_{N_0}\sigma^k\,.
\]
As $\sigma>1$ and $P_{N_0}>0$, the sequence $P_{N_0+k}\geq P_{N_0}\sigma^k$ diverges to infinity as $k\to \infty$. This contradicts the fact that $X$ is compact invariant set (specifically $P(t)\leq C_0~~\forall t\geq 0$) and also contradicts that $x_n\to E_0$. Thus, $f$ is weakly uniform persistent with respect to $(X_0,\partial X_0)$. Since $f$ has a global attractor and is weakly uniform persistent, Theorem 1.3.3 from \cite{zhao} implies that $f$ is uniformly persistent, \emph{i.\,e.}, there exists $\epsilon>0$ such that
\begin{equation*}
\liminf_{n\to \infty}{P_{x}(nT_p)}\geq \epsilon~~\forall x\in X_0\,.
\end{equation*}

Finally, we know that the maximum per-capita death rate is finite. Letting it be $L$ we obtain
\begin{equation*}
\frac{\dot{P}}{P}\geq -L\,.
\end{equation*}
Let $n\in \mathbb{N}$, the above implies that for all $t\in\left[nT_p,nT_p+T_p\right]$ we get
\begin{equation*}
P(t)\geq P(nT_p)\exp(-L(t-nT_p))\geq P(nT_p)\exp(-LT_p)\geq \epsilon \exp(-LT_p)=\eta >0\,.
\end{equation*}
The above holds true every year, therefore, we see that persistence is also uniform in the continuous system, \emph{i.\,e.},
\begin{equation*}
\liminf_{t\to \infty}{P_{x}(t)}\geq \eta~~\forall x\in X_0\,.
\end{equation*}
\qed

\section*{Acknowledgments}

Work funded by the Franco-Chilean Binational Center of Artificial Intelligence, ANID Strengthening R\&D capabilities Program CTI230007 Inria Chile, and Inria Challenge OcéanIA (desc. num 14500). Pablo A. Marquet acknowledges support from Centro de Modelamiento Matem\'atico (CMM), Grant FB210005, BASAL funds for centers of excellence from ANID-Chile, and Proyecto Exploraci\'on 13220168. The authors are also grateful to Aurora Llanos and Ignacio Díaz for their valuable support during the final stages of this project; their assistance with technical adjustments and their careful proofreading were instrumental in the preparation of this manuscript.

\section*{Author Contributions}

All authors have made substantial intellectual contributions to the conception, execution, or design of the work. All authors have read and approved the final manuscript. Contributions, as specified by the Contributor Role Taxonomy (CRediT) can be summarized as:
\begin{itemize}
    \item Conceptualization: Matías Neto, Pablo A. Marquet, Mara Freilich, Luis Martí, Nayat Sanchez-Pi.
    \item Data curation: Mara Freilich, Matías Neto, Luis Martí.
    \item Formal Analysis: Matias Neto.
    \item Funding acquisition: Nayat Sanchez-Pi.
    \item Investigation:  Matías Neto, Pablo A. Marquet, Mara Freilich, Luis Martí, Nayat Sanchez-Pi.
    \item Methodology: Matías Neto, Pablo A. Marquet, Mara Freilich.
    \item Project administration: Nayat Sanchez-Pi.
    \item Resources: Matías Neto, Pablo A. Marquet, Mara Freilich, Luis Martí, Nayat Sanchez-Pi.
    \item Software: Matías Neto, Luis Martí.
    \item Supervision: Pablo A. Marquet, Nayat Sanchez-Pi.
    \item Validation: Matías Neto, Pablo A. Marquet, Mara Freilich, Luis Martí, Nayat Sanchez-Pi.
    \item Visualization: Matías Neto, Luis Martí.
    \item Writing—Original Draft Preparation: Matías Neto.
    \item Writing—Critical Review and Editing: Pablo A. Marquet, Mara Freilich, Luis Martí, Nayat Sanchez-Pi.
\end{itemize}
% All authors agree to be responsible for all aspects of the work, ensuring that any questions related to accuracy or integrity are investigated and resolved appropriately. A large language model (LLM) based on artificial intelligence was used to edit and polish the written text in order to improve grammar and overall style.
% \subsection{Conflicts of Interest}
% The authors declare there are no conflicts of interest
% \subsection{Ethics Statement}
% The authors confirm that this study is based on theoretical mathematical modeling and numerical simulations using publicly available climate model data. Therefore, the work did not involve human participants, animals, or primary data collection that would require review by an institutional review board or ethics committee.

\section*{Access to code} 
Code is available at \url{https://github.com/Inria-Chile/NPZ-T_Persistence}.

\bibliographystyle{siamplain}
\bibliography{npz-t}

\end{document}